\documentclass{amsart}
\usepackage{amsmath,amsfonts,amscd,amsthm,graphicx,amssymb,epic,eepic}

\theoremstyle{plain}
\newtheorem{theorem}{Theorem}
\newtheorem{proposition}{Proposition}
\newtheorem{cor}{Corollary}
\newtheorem{lemma}{Lemma}

\numberwithin{equation}{section}

\newcommand{\eps}{\varepsilon}

\newcommand{\Z}{\mathbb Z}
\newcommand{\N}{\mathbb N}
\newcommand{\R}{\mathbb R}

\renewcommand{\P}{\mathbb P}

\newcommand{\EE}{\mathcal E}

\renewcommand{\H}{\mathbb H}

\newcommand{\NN}{\mathcal N}
\newcommand{\MM}{\mathcal M}

\title[Angles between geodesic rays associated with
hyperbolic lattice points]{On the distribution of angles between
geodesic rays \\ associated with hyperbolic lattice points}
\author{Florin P. Boca}

\address{Department of Mathematics, University of Illinois, 1409 W. Green Street, Urbana, IL 61801, USA}

\address{Institute of Mathematics ``Simion Stoilow" of the Romanian
Academy, P.O. Box 1-764, RO-014700 Bucharest, Romania}

\email{fboca@math.uiuc.edu}
\subjclass[2000]{Primary: 11P21; Secondary: 11L05, 30F35, 51M09.}

\date{January 1, 2007}

\thanks{}

\begin{document}

%\maketitle

\begin{abstract}
For every two points $z_0,z_1$ in the upper-half plane $\H$,
consider all elements $\gamma$ in the principal congruence group
$\Gamma(N)$, acting on $\H$ by fractional linear transformations,
such that the hyperbolic distance between $z_1$ and $\gamma z_0$
is at most $R>0$. We study the distribution of angles between the
geodesic rays $[z_1,\gamma z_0]$ as $R\rightarrow\infty$, proving
that the limiting distribution exists independently of $N$ and
explicitly computing it. When $z_1=z_0$ this is found to be the
uniform distribution on the interval
$\big[-\frac{\pi}{2},\frac{\pi}{2}\big]$.
\end{abstract}

\maketitle

\section{Introduction}
In this paper the group $SL_2(\R)$ acts on the upper half-plane
$\H$ by linear fractional transformations $z\mapsto g
z=\frac{az+b}{cz+d}$, $g=\left(\begin{smallmatrix} a & b \\ c & d
\end{smallmatrix}\right)\in SL_2(\R)$, $z\in\H$.
The hyperbolic ball $B(z_0,R)=\{ z\in\H:\varrho (z_0,z)\leq R\}$
of center $z_0=x_0+iy_0\in\H$ and radius $R$ coincides with the
Euclidean ball of center $x_0+iy_0\cosh R$ and radius $y_0\sinh R
\sim \frac{1}{2}y_0e^R$. Let $\Gamma$ be a discrete subgroup of
$SL_2(\R)$. The \emph{hyperbolic circle problem} of estimating for
fixed $z_0,z_1\in\H$ and $R\rightarrow \infty$ the cardinality of
the set $\Gamma_{z_0,R}=\{ \gamma \in\Gamma:\gamma z_0\in
B(z_0,R)\}$, or slightly more generally of $\{ \gamma\in\Gamma:
\varrho (\gamma z_0,z_1)\leq R\}$, has been thoroughly studied
with various methods (see, e.g., \cite{De,Go,Hej,HZ,Iw,LP,Pa}, and
\cite{LP,Ma,MaSh} for some higher dimensional analogs of the
problem).

We consider another natural problem concerning the distribution of
hyperbolic lattice points in angular sectors. For $z_0,z_1\in \H$
and $g\in SL_2(\R)$, let $\theta_{z_0,z_1}(g)\in
[-\frac{\pi}{2},\frac{\pi}{2}]$ denote the angle between the
geodesic ray $[z_1,g z_0]$ and the vertical geodesic
$[z_1,\infty]$. Given a compact set $\Omega \subset \H$ and a
number $\omega\in [-\frac{\pi}{2},\frac{\pi}{2}]$, the proportion
of points in the $\Gamma$-orbit of $z_0$ inside $\Omega$ such that
$\theta_{z_0,z_1}(\gamma)\leq\omega$ is given by
\begin{equation*}
\P_{\Gamma,\Omega,z_0,z_1}(\omega)=\frac{\# \{ \gamma\in\Gamma :
\gamma z_0\in\Omega,\ \theta_{z_0,z_1}(\gamma)\leq \omega \}}{\#
\{ \gamma\in\Gamma :\gamma z_0\in\Omega\}}.
\end{equation*}
It is natural to investigate the existence of the limiting
distribution
\begin{equation*}
\P_{\Gamma,z_0,z_1}(\omega)=\lim\limits_{R\rightarrow\infty}
\P_{\Gamma,B(z_0,R),z_0,z_1} (\omega)=\lim\limits_{R\rightarrow
\infty} \frac{\# \{ \gamma\in \Gamma_{z_0,R}: \theta_{z_0,z_1}
(\gamma)\leq\omega\}}{\# \Gamma_{z_0,R}}.
\end{equation*}

In this paper we consider the case where
\begin{equation*}
\Gamma=\Gamma(N)=\left\{ \left( \begin{matrix} a & b \\ c & d
\end{matrix}\right)\in SL_2(\Z): a,d\equiv 1,\ b,c\equiv
0\hspace{-6pt}\pmod{N} \right\}
\end{equation*}
is the \emph{principal congruence subgroup of level $N$}, which is
the kernel of the natural surjective morphism $SL_2(\Z)\rightarrow
SL_2(\Z_N)$. This is a normal subgroup of $\Gamma(1)=SL_2(\Z)$ of
index
\begin{equation}\label{1.1}
[\Gamma(1):\Gamma(N)]=N^3 \prod\limits_{\substack{p\mid N \\ p\
\mathrm{prime}}} (1-p^{-2}).
\end{equation}
For every $g=\left(\begin{smallmatrix} A & B \\ C & D
\end{smallmatrix}\right)\in SL_2 (\R)$ the hyperbolic distance $\varrho (i,g i)$ is given by
\begin{equation}\label{1.2}
\cosh \varrho (i,g i)=1+\frac{\vert i-g
i\vert^2}{2\operatorname{Im} (gi)}=\frac{A^2+B^2+C^2+D^2}{2}.
\end{equation}
Denote
\begin{equation}\label{1.3}
C_N=\sum\limits_{\substack{n\geq 1 \\
(n,N)=1}} \frac{\mu (n)}{n^2}=\prod\limits_{p\nmid N}
(1-p^{-2})=\frac{1}{\zeta(2)} \prod\limits_{p\mid
N}(1-p^{-2})^{-1}.
\end{equation}

For every $z_0=x_0+iy_0, z_1=x_1+iy_1\in\H$, denote
$x_*=\frac{x_1-x_0}{y_0}$, $y_*=\frac{y_1}{y_0}$, and consider the
continuous function $\Xi_{x_*,y_*}$ on
$[-\frac{\pi}{2},\frac{\pi}{2}]$ defined by
\begin{equation}\label{1.4}
\begin{split}
\Xi_{x_*,y_*}(\omega) = & \frac{1}{\pi}\arctan \left(
x_*+y_*\tan\frac{\omega}{2}\right) +\frac{1}{\pi}\arctan \left(
x_*-y_*\cot\frac{\omega}{2}\right) \\ & -\frac{1}{\pi}\arctan
(x_*+y_*)-\frac{1}{\pi}\arctan (x_*-y_*)+\begin{cases} 1 &
\mbox{\rm if $\omega >0,$} \\ 0 & \mbox{\rm if $\omega <0.$}
\end{cases}
\end{split}
\end{equation}

The main result of this paper is

\begin{theorem}\label{T1}
For every positive integer $N$ and $z_0=x_0+iy_0,
z_1=x_1+iy_1\in\H$, as $R\rightarrow \infty$,
\begin{equation}\label{1.5}
\# \left\{ \gamma\in \Gamma (N)_{z_0,R}: -\frac{\pi}{2}\leq
\theta_{z_0,z_1}(\gamma)\leq\omega\right\} = \frac{\pi^2 C_N
\Xi_{x_*,y_*}(\omega)}{N^3}e^R +O_{\eps,N,z_0,z_1} \left(
e^{(7/8+\eps)R}\right).
\end{equation}
In particular the limiting distribution $\P_{\Gamma(N),z_0,z_1}$
exists and is given by
\begin{equation*}
\P_{\Gamma(N),z_0,z_1}(\omega)=\frac{1}{\pi}\int_{-\pi/2}^\omega
\varrho_{z_0,z_1}(t)\, dt,\qquad \omega\in
\left[-\frac{\pi}{2},\frac{\pi}{2}\right],
\end{equation*}
where
\begin{equation*}
\varrho_{z_0,z_1}(t)=\frac{2y_0y_1\big(y_0^2+(x_1-x_0)^2+y_1^2\big)}{\big(
y_0^2+(x_1-x_0)^2+y_1^2\big)^2-\big( (y_0^2+(x_1-x_0)^2-y_1^2)\cos
t+2(x_1-x_0)y_1 \sin t\big)^2}.
\end{equation*}
\end{theorem}

Taking $z_1=z_0$ we infer

\begin{cor}\label{C1}
The angles $\theta_{z_0,z_0}(\gamma)$,
$\gamma\in\Gamma(N)_{z_0,R}$, are uniformly distributed as
$R\rightarrow\infty$.
\end{cor}

The converse is also seen to be true, so that the angles
$\theta_{z_0,z_1}(\gamma)$ are uniformly distributed as
$R\rightarrow\infty$ if and only if $z_1=z_0$. In the Euclidean
situation these angles are uniformly distributed regardless of the
choice of $z_1$ and $z_0$.

Our method of proof is number theoretical and relies on the Weil
bound for Kloosterman sums \cite{We}, as previously used (for
instance) in \cite{Bo,BCZ,BZ,Fu,HZ}. In the process we also
derive, as a consequence of the proof of Theorem \ref{T1}, an
asymptotic formula for the number of hyperbolic lattice points in
large balls.

\begin{cor}\label{C2}
For every positive integer $N$ and every $z_0\in\H$, as
$R\rightarrow\infty$,
\begin{equation}\label{1.6}
\# \Gamma(N)_{z_0,R}=\frac{6e^R}{[\Gamma (1):\Gamma
(N)]}+O_{\eps,N,z_0}\left(e^{( 7/8+\eps)R}\right).
\end{equation}
\end{cor}

Denoting by $\mu$ the hyperbolic area in $\H$, the main term in
\eqref{1.5} is $\sim \frac{2\mu(B(z_0,R))}{\mu
(\Gamma(N)\backslash\H)}$ as $R\rightarrow \infty$.

For $N=1$ formula \eqref{1.6} has been proved using Kloosterman
sum estimates in \cite{HZ}. Better error terms with exponent as
low as $\frac{2}{3}$ can be obtained using Selberg's theory on the
spectral decomposition of $L^2(\Gamma(N)\backslash\H)$ (see
\cite{Pa} for exponent $\frac{3}{4}$ and \cite{Iw} for exponent
$\frac{2}{3}$) and lower bounds for the first eigenvalue of the
Laplacian on $\Gamma(N)\backslash\H$ (see \cite{Sel}, \cite{Iw},
and \cite{Sha} for a review of recent developments). Similar
results hold when $\Gamma(N)$ is replaced by any of the congruence
groups $\Gamma_0(N)=\{\gamma\in\Gamma(1): c\equiv
0\hspace{-3pt}\pmod{N}\}$ or $\Gamma_1(N)=\{ \gamma\in \Gamma(1):
a,d\equiv 1, c\equiv 0\hspace{-3pt}\pmod{N}\}$, or when
$\Gamma(N)_{z_0,R}$ is replaced by $\{
\gamma\in\Gamma(N):\varrho(\gamma z_0,z_1)\leq R\}$ for fixed
$z_0,z_1\in \H$.

There are two natural problems that arise in this context. It
would be interesting to know how large is the class of discrete
subgroups of $SL_2(\R)$ for which the analogue of Theorem \ref{T1}
holds. It would also be interesting to study the spacing
statistics (both consecutive spacings and correlations) of these
angles when $z_0=z_1$.

\section{Reducing the problem to a counting problem}
Given $z_0=x_0+iy_0\in\H$ and $\gamma=\left(\begin{smallmatrix} a & b \\
c & d \end{smallmatrix}\right)\in\Gamma(1)$, consider
\begin{equation*}
g_0=\left( \begin{matrix} \sqrt{y_0} & \frac{x_0}{\sqrt{y_0}} \\
0 & \frac{1}{\sqrt{y_0}} \end{matrix} \right),\quad g =g_0^{-1}
\gamma g_0=\left( \begin{matrix} A & B \\ C & D
\end{matrix}\right)\ \in\ SL_2(\R),
\end{equation*}
with
\begin{equation}\label{2.1}
A=a-cx_0,\quad B=\frac{(a-cx_0)x_0+b-dx_0}{y_0},\quad C=cy_0,\quad
D=cx_0+d.
\end{equation}

Since $g_0 i=z_0$ we have
\begin{equation}\label{2.2}
\cosh \varrho(z_0,\gamma z_0)=\cosh \varrho (g_0 i,g_0 g i)=\cosh
\varrho (i,g i)=\frac{A^2+B^2+C^2+D^2}{2}.
\end{equation}

Take $Q^2=2\cosh R\sim e^R$. As a result of \eqref{2.2} we are
interested in those $\gamma\in\Gamma(N)$ for which
$A^2+B^2+C^2+D^2\leq Q^2$. The only matrices $\gamma\in \Gamma
(1)$ with $c=0$ are $\pm I_2$ and as a result we can assume next
that $C\neq 0$. We will also assume that $A\neq 0$.

The geodesic joining the points $z_*=x_*+iy_*\in\H$ and $gi$,
$g=\left(\begin{smallmatrix} A & B \\ C & D
\end{smallmatrix}\right)\in SL_2(\R)$, is the half-circle of
center $\alpha$ and radius $r$, where
\begin{equation*}
\vert\alpha-z_*\vert=\vert \alpha -gi\vert=r.
\end{equation*}
This gives
\begin{equation*}
\vert \alpha-x_*-iy_* \vert^2 =\left|
\alpha-\frac{iA+B}{iC+D}\right|^2 = \frac{\vert
i(C\alpha-A)+D\alpha -B\vert^2}{\vert iC+D\vert^2},
\end{equation*}
and after cancelling out the terms containing $\alpha^2$ we obtain
\begin{equation*}
2\alpha (x_* E-F)=(x_*^2+y_*^2)E-G,
\end{equation*}
with
\begin{equation*}
E=C^2+D^2,\qquad F=AC+BD,\qquad G=A^2+B^2,
\end{equation*}
leading to
\begin{equation*}
\tan\theta_{i,z_*}(g)=\frac{y_*}{x_*-\alpha}=\frac{2y_*(F-x_*
E)}{(y_*^2-x_*^2)E+2x_* F-G}.
\end{equation*}
We will keep $z_0$ and $z_1$ fixed throughout. Taking
\begin{equation*}
z_*=g_0^{-1} z_1=\frac{x_1-x_0+iy_1}{y_0}
\end{equation*}
we have $g_0(x_*+it)=x_1+iy_0 t$, $t>0$, so that
\begin{equation*}
\theta_{z_0,z_1}(\gamma)=\measuredangle [x_1+i\infty,z_1,\gamma
z_0]=\measuredangle [g_0(x_*+i\infty), g_0 z_*,g_0 g
i]=\measuredangle [x_*+i\infty,z_*,gi]=\theta_{i,z_*}(g),
\end{equation*}
and therefore
\begin{equation}\label{2.3}
\tan\theta_{z_0,z_1}(\gamma)=\frac{y_*}{x_*-\alpha}=\frac{2y_*(F-x_*
E)}{(y_*^2-x_*^2)E+2x_* F-G}.
\end{equation}

When $\vert A\vert\leq\vert D\vert$ we use
\begin{equation*}
\begin{split}
\left|  \frac{F}{E}-\frac{A}{C}\right| &
=\frac{\vert D\vert}{\vert C\vert(C^2+D^2)} \leq \frac{1}{2C^2}, \\
\left| \frac{G}{E}-\frac{A^2}{C^2}\right| & =\frac{\vert
BC+AD\vert}{C^2(C^2+D^2)} =\frac{\vert 2AD-1\vert}{C^2(C^2+D^2)}
\leq \frac{2}{C^2}+\frac{1}{C^4},
\end{split}
\end{equation*}
to derive
\begin{equation}\label{2.4}
\tan\theta_{z_0,z_1}(\gamma)
=\frac{2y_*\big(\frac{F}{E}-x_*\big)}{y_*^2-x_*^2+2x_*\frac{F}{E}-\frac{G}{E}}
=\frac{2y_*\big(\frac{A}{C}-x_*\big)
+O_{z_*}\big(\frac{1}{C^2}\big)}{y_*^2-\big(\frac{A}{C}-x_*\big)^2
+O_{z_*}\big(\frac{1}{C^2}+\frac{1}{C^4}\big)}.
\end{equation}

When $\vert D\vert\leq\vert A\vert$ we use
\begin{equation*}
\begin{split}
\left| \frac{F}{G}-\frac{C}{A}\right| &
=\frac{\vert B\vert}{\vert A\vert(A^2+B^2)}\leq\frac{1}{2A^2},\\
\left| \frac{E}{G}-\frac{C^2}{A^2}\right| & =\frac{\vert
2AD-1\vert}{A^2(A^2+B^2)} \leq \frac{2}{A^2}+\frac{1}{A^4},
\end{split}
\end{equation*}
to derive
\begin{equation}\label{2.5}
\begin{split}
\tan\theta_{z_0,z_1}(\gamma) &
=\frac{2y_*\big(\frac{F}{G}-x_*\frac{E}{G}\big)}{(y_*^2-x_*^2)
\frac{E}{G}+2x_*\frac{F}{G}-1}
=\frac{2y_*\big(\frac{C}{A}-x_*\frac{C^2}{A^2}\big)
+O_{z_*}\big(\frac{1}{A^2}+\frac{1}{A^4}\big)}{(y_*^2-x_*^2)\frac{C^2}{A^2}
+2x_*\frac{C}{A}-1+O_{z_*}\big(\frac{1}{A^2}+\frac{1}{A^4}\big)} \\
&
=\frac{2y_*\big(\frac{A}{C}-x_*\big)+O_{z_*}\big(\frac{1}{C^2}+\frac{1}{A^2
C^2}\big)}{y_*^2-\big(\frac{A}{C}-x_*\big)^2+O_{z_*}
\big(\frac{1}{C^2}+\frac{1}{A^2C^2}\big)}.
\end{split}
\end{equation}

For $\lambda >0$ set
\begin{equation*}
-\alpha_1:=\frac{-1-\sqrt{1+\lambda^2}}{\lambda}
<-1<0<\alpha_2:=\frac{-1+\sqrt{1+\lambda^2}}{\lambda}=\frac{1}{\alpha_1}<1.
\end{equation*}
For $\lambda <0$ set
\begin{equation*}
-1<\alpha_1^*:=\frac{1-\sqrt{1+\lambda^2}}{\vert\lambda\vert}<0<1<\alpha_2^*:=
\frac{1+\sqrt{1+\lambda^2}}{\vert\lambda\vert}=-\frac{1}{\alpha_1^*}.
\end{equation*}
Letting $\lambda=\tan\omega$,
$\omega\in\big(-\frac{\pi}{2},\frac{\pi}{2}\big)$, we have
$\alpha_1=\cot\frac{\omega}{2}$, $\alpha_2=\tan\frac{\omega}{2}$
for $\omega>0$, and $\alpha_1^*=\tan\frac{\omega}{2}$,
$\alpha_2^*=-\cot\frac{\omega}{2}$ for $\omega <0$. A plain
calculation gives
\begin{equation}\label{2.6}
\frac{2y_*(X-x_*)}{y_*^2-(X-x_*)^2}<\lambda \ \Longleftrightarrow\
X-x_*\in\mathfrak{S}(y_*,\lambda),
\end{equation}
with
\begin{equation}\label{2.7}
\mathfrak{S}(y_*,\lambda)=\begin{cases} (-\infty,-y_*\alpha_1
)\cup (-y_*,y_*\alpha_2 )\cup (y_*,\infty) &
\mbox{\rm if $\lambda >0$,} \\
(-y_*,0)\cup (y_*,\infty) & \mbox{\rm if $\lambda=0$,} \\
(-y_*,y_* \alpha_1^*)\cup (y_*,y_*\alpha_2^*) & \mbox{\rm if
$\lambda <0$.}
\end{cases}
\end{equation}

For fixed $\lambda\in\R$, $z_*\in\H$, and $\vert\eps_1\vert$,
$\vert\eps_2\vert$ small, the roots $X_\pm (\eps_2)$ of
$y_*^2-(X-x_*)^2+\eps_2=0$ and $\widetilde{X}_\pm (\eps_1,\eps_2)$
of $2y_*(X-x_*)+\eps_1-\lambda\big( y_*^2-(X-x_*)^2+\eps_2\big)=0$
satisfy
\begin{equation*}
\vert X_\pm (\eps_2)-X_\pm (0)\vert=\Big|
\sqrt{y_*^2+\eps_2}-y_*\Big| \leq \frac{\vert\eps_2\vert}{y_*},
\end{equation*}
and respectively
\begin{equation*}
\left| \widetilde{X}_\pm (\eps_1,\eps_2)-\widetilde{X}_\pm
(0,0)\right| =\frac{\vert \eps_1-\lambda
\eps_2\vert}{y_*\sqrt{1+\lambda^2}+\sqrt{y_*^2
(1+\lambda^2)-\lambda\eps_1+\lambda^2\eps_2}} \leq \frac{\vert
\eps_1 -\lambda \eps_2\vert}{y_*\sqrt{1+\lambda^2}}\leq
\frac{\sqrt{\eps_1^2+\eps_2^2}}{y_*}.
\end{equation*}
In conjunction with \eqref{2.4}--\eqref{2.7} this shows, in both
cases $\vert A\vert\leq \vert D\vert$ and $\vert D\vert\leq \vert
A\vert$, that there is a constant $K_1=K_1(z_*) >0$ such that, for
any $\gamma\in\Gamma (N)$,
\begin{equation}\label{2.8}
\tan \theta_{z_0,z_1}(\gamma) \leq\lambda \ \Longrightarrow\
\frac{A}{C}\in
x_*+\mathfrak{S}(y_*,\lambda)+\big[-K_1H(\gamma),K_1H(\gamma)\big],
\end{equation}
where $H(\gamma)=\frac{1}{C^2}+\frac{1}{A^2C^2}+\frac{1}{C^4}$.

We wish to discard those $\gamma$ for which one of $\vert A\vert$,
$\vert B\vert$, $\vert C\vert$, $\vert D\vert$ is small. Note
first that, as a result of \eqref{2.1}, there is a constant
$K_0=K_0(z_0)$ such that $a^2+b^2+c^2+d^2\leq K_0 Q^2$ whenever
$A^2+B^2+C^2+D^2\leq Q^2$. For every $K>0$ let
$\EE_A(K)=\EE_{A,Q,z_0}(K)$ denote the number of
$\gamma\in\Gamma(1)$ for which $A^2+B^2+C^2+D^2\leq Q^2$ and
$\vert A\vert=\vert a-cx_0\vert \leq K$. Define similarly
$\EE_B(K)$, $\EE_C (K)$, $\EE_D(K)$.

\begin{lemma}\label{L1}
{\em (i)} For every $z_0\in\H$ and $K\geq 1$
\begin{equation*}
\max \left\{ \EE_A(K),\EE_C(K),\EE_D(K)\right\} \ll_{z_0} KQ\log Q
\qquad (Q\rightarrow \infty).
\end{equation*}

{\em (ii)}. For every $z_0\in\H$ and $\alpha\in (0,1)$
\begin{equation*}
\EE_B(Q^\alpha)\ll_{z_0} Q^{(3+\alpha)/2}\log Q \qquad
(Q\rightarrow\infty).
\end{equation*}
\end{lemma}

\begin{proof}
(i) The congruence $bc=1\pmod{\vert a\vert}$ shows, for fixed $c$
and $a\neq 0$, that the integer $b$ is uniquely determined
$\pmod{\vert a\vert}$, so it takes $\ll \frac{Q}{\vert a\vert}$
values. This gives
\begin{equation*}
\EE_C(K)\ll 2+\left(\frac{2K}{y_0}+1\right)\sum\limits_{1\leq
\vert a\vert \leq K_0 Q} \frac{Q}{\vert a\vert} \ll_{z_0} KQ\log
Q.
\end{equation*}

To prove $\EE_A(K)\ll_{z_0} KQ\log Q$ note that, for fixed $c\in
[-K_0 Q,K_0 Q]$, there are at most $2K+1$ integers $a$ such that
$\vert a-cx_0\vert\leq K$. For each such $a$, the congruence
$ad=1\pmod{\vert c\vert}$ uniquely determines $d\pmod{\vert
c\vert}$, so the number of admissible triples $(a,d,b)$ is $\ll
\frac{KQ}{\vert c\vert}$, and summing over $c$ we find as above
$\EE_A(K)\ll_{z_0} KQ\log Q$. The proof of $\EE_D(K)\ll_{z_0}
KQ\log Q$ is similar.

(ii) Let $\EE_A(K)^c$, respectively $\EE_D(K)^c$, denote the
complement of $\EE_A(K)$, respectively $\EE_D(K)$, in $\{
\gamma\in \Gamma(1): A^2+B^2+C^2+D^2\leq Q^2\}$. Write
$\alpha=2\alpha^\prime -1$, $\frac{1}{2}<\alpha^\prime<1$, so that
$1+\alpha^\prime=\frac{3+\alpha}{2}$. For every
$\gamma\in\EE_A(Q^{\alpha^\prime}+1)^c\cap
\EE_D(Q^{\alpha^\prime}+1)^c$ we have
\begin{equation*}
\vert B\vert=\frac{\vert AD-1\vert}{\vert C\vert}
>\frac{Q^{2\alpha^\prime}}{Q}=Q^{\alpha},
\end{equation*}
showing that $\EE_B(Q^{\alpha}) \subseteq
\EE_A(Q^{\alpha^\prime}+1)\cup\EE_D(Q^{\alpha^\prime}+1)$, and so
$\EE_D(Q^\alpha)\ll_{z_0} Q^{1+\alpha^\prime}\log Q$.
\end{proof}

Note also that
\begin{equation}\label{2.9}
\begin{split}
\left| (A^2+B^2+C^2+D^2)-(C^2+A^2)\bigg( 1
+\frac{D^2}{C^2}\bigg)\right| & =\frac{\vert AD+BC\vert}{C^2}
=\frac{\vert 2BC+1\vert}{C^2} \\ & \leq \frac{2\vert B\vert}{\vert
C\vert}+\frac{1}{C^2} \ll_{z_0} \frac{Q}{\vert
c\vert}+\frac{1}{c^2}\ll Q.
\end{split}
\end{equation}

The relations \eqref{2.8} and \eqref{2.9} lead us to estimate the
number
\begin{equation}\label{2.10}
\mathfrak{N}_Q (N,z_0;\beta):=\#\left\{ \gamma\in\Gamma(N):
\frac{A}{C}\leq\beta,\ (C^2+A^2)\left(1+\frac{D^2}{C^2}\right)\leq
Q^2\right\}\qquad (Q\rightarrow\infty).
\end{equation}

\section{Some counting in $\Gamma (N)$}
In this section we prove some counting results which will be
further used in the proof of Theorem \ref{T1} in the next section.
Let $c$ and $N\geq 1$ be integers and consider the sum
\begin{equation*}
\Phi_N (c): =\sum\limits_{\substack{n\mid c \\ (n,N)=1}} \frac{\mu
(n)}{n}.
\end{equation*}

We first estimate the number
\begin{equation*}
\NN_{c,N}(I_1\times I_2):=\# \big\{ (a,d)\in I_1 \times I_2 :
a\equiv 1, \ d\equiv 1\hspace{-6pt} \pmod{N},\ ad\equiv 1
\hspace{-6pt} \pmod{Nc}\big\},
\end{equation*}
with fixed $N$ and $c$, and with $a$ and $d$ in prescribed (short)
intervals. The next result extends Lemma 1.6 in \cite{BCZ} from
$\Gamma(1)$ to $\Gamma(N)$.

\begin{proposition}\label{P1}
For a fixed positive integer $N$ and intervals $I_1,I_2$ of length
less than $\vert c\vert$
\begin{equation*}
\NN_{c,N}(I_1\times I_2)=\frac{\Phi_N(c)}{\vert c\vert N^2}
\,\vert I_1\vert\,\vert I_2\vert+O_{\eps,N} ( \vert
c\vert^{1/2+\eps} )\qquad (\vert c\vert \rightarrow\infty).
\end{equation*}
\end{proposition}

\begin{proof}
Replacing $(b,c)$ by $(-b,-c)$ we can assume $c>0$. In this case
we write
\begin{equation*}
\NN_{c,N} (I_1\times I_2) =\frac{1}{N
c}\sum\limits_{\substack{x\in I_1\\ (x,Nc)=1 \\
x\equiv 1\hspace{-6pt} \pmod{N}}} \sum\limits_{\substack{y\in I_2
\\ y\equiv 1\hspace{-6pt}\pmod{N}}}
\sum\limits_{k\hspace{-6pt} \pmod{Nc}} e\left(
\frac{k(y-\bar{x})}{Nc}\right) =\MM +\EE,
\end{equation*}
where $\bar{x}$ is the multiplicative inverse of
$x\hspace{-3pt}\pmod{Nc}$ and $e(t)=\exp (2\pi it)$. The
contribution
\begin{equation}\label{3.1}
\MM=\frac{1}{N c} \sum\limits_{\substack{x\in I_1\\ (x,Nc)=1 \\
x\equiv 1\hspace{-6pt} \pmod{N}}}
\sum\limits_{\substack{y\in I_2 \\
y\equiv 1\hspace{-6pt} \pmod{N}}} \sum\limits_{0\leq \ell<N}
e\left( \frac{\ell (y-\bar{x})}{N}\right)
\end{equation}
of terms with $c\mid k$ to $\NN_{c,N}(I_1,I_2)$ will be treated as
a main term, while the contribution
\begin{equation}\label{3.2}
\EE=\frac{1}{N c} \sum\limits_{\substack{0\leq k<Nc \\
c\nmid k}}\sum\limits_{\substack{y\in I_2\\ y\equiv 1\hspace{-6pt}
\pmod{N}}} e\left(\frac{ky}{Nc}\right)
\sum\limits_{\substack{x\in I_1 \\ (x,Nc)=1 \\
x\equiv 1\hspace{-6pt} \pmod{N}}} e\left(
-\frac{k\bar{x}}{Nc}\right)
\end{equation}
of terms with $c\nmid k$ will be treated as an error term.

To estimate $\EE$ consider for $I$ interval and $q\in \N$,
$m,n\in\Z$, the incomplete Kloosterman sum
\begin{equation*}
S_I (m,n;q):=\sum\limits_{\substack{a\in I\\ (a,q)=1}} e\left(
\frac{ma+n\bar{a}}{q}\right),
\end{equation*}
where $\bar{a}$ is the multiplicative inverse of $a\pmod{q}$. The
complete Kloosterman sum $S(m,n;q)$ is just $S_{[0,q-1]}(m,n;q)$.
For any interval $I\subseteq [0,q-1]$ and integers $m,n$, not both
divisible by $q$, the Weil bound on Kloosterman sums leads (cf.,
e.g., \cite[Lemma 1.6]{BCZ}) to
\begin{equation}\label{3.3}
\vert S_I (m,n;q)\vert \ll_\eps (n,q)^{1/2} q^{1/2+\eps}.
\end{equation}

Writing now the inner sum in \eqref{3.2} as
\begin{equation*}
\sum\limits_{\substack{x\in I_1 \\ (x,Nc)=1}} e\left(
-\frac{k\bar{x}}{Nc}\right) \frac{1}{N}
\sum\limits_{s\hspace{-6pt} \pmod{N}} e\left(
\frac{s(x-1)}{N}\right) =\frac{1}{N}\sum\limits_{s\hspace{-6pt}
\pmod{N}} e\left( -\frac{s}{N}\right) S_{I_1} (cs,-k;Nc)
\end{equation*}
and applying \eqref{3.3} we find
\begin{equation*}
\vert \EE \vert \ll_\eps \frac{(Nc)^{1/2+\eps}}{N c}
\sum\limits_{\substack{0\leq k<Nc \\ c\nmid k}} (k,Nc)^{1/2}\Bigg|
\sum\limits_{\substack{y\in I_2
\\ y\equiv 1\hspace{-6pt} \pmod{N}}} e\bigg( \frac{ky}{Nc}\bigg) \Bigg| .
\end{equation*}
Treating the inner sum above as a geometric progression of ratio
$e( \frac{k}{c})$ and using the inequality $\vert\sin \pi
t\vert\geq 2\| t\|=2\operatorname{dist}(t,\Z)$, $t\in\R$, the
inner sum above is $\leq \min\big\{ \vert I_2\vert,\frac{1}{2\|
k/c\|}\big\}$. Employing also the inequality $(k,Nc)\leq (k,c)N$
we further find
\begin{equation*}
\begin{split}
\vert\EE\vert & \ll_\eps N^{1+\eps} \frac{c^{1/2+\eps}}{c}
\sum\limits_{0<\ell<c} \frac{(\ell,c)^{1/2}}{2\big\|
\frac{\ell}{c}\big\|} \ll N^{1+\eps} c^{-1/2+\eps}
\sum\limits_{d\mid c} \sum\limits_{m\leq\frac{c}{2d}}
d^{1/2}\,\frac{c}{dm} \\
& \leq N^{1+\eps} c^{1/2+\eps} \sum\limits_{d\mid c} d^{-1/2} \log
c \ll_{\eps,N} c^{1/2+2\eps}.
\end{split}
\end{equation*}

Concerning the main term $\MM$, from $x\bar{x}=1\hspace{-3pt}
\pmod{N}$ and $x=1\hspace{-3pt} \pmod{N}$ we infer
$\bar{x}=1\hspace{-3pt} \pmod{N}$, and so $N\mid (y-\bar{x})$. The
inner sum in \eqref{3.1} is equal to $N$ and we get
\begin{equation*}
\MM=\frac{1}{c}\sum\limits_{\substack{x\in I_1 \\ (x,Nc)=1 \\
x\equiv 1\hspace{-6pt}\pmod{N}}}1 \sum\limits_{\substack{y\in I_2 \\
y\equiv 1\hspace{-6pt}\pmod{N}}} 1=\frac{1}{c}\left( \frac{\vert
I_2\vert}{N}+O(1)\right)\sum\limits_{\substack{x\in I_1 \\ (x,Nc)=1 \\
x\equiv 1\hspace{-6pt} \pmod{N}}} 1.
\end{equation*}
Using $x=1\hspace{-3pt} \pmod{N}$ and M\" obius summation, the
latter sum above can also be expressed as
\begin{equation*}
\begin{split}
\sum\limits_{\substack{ x\in I_1 \\ x\equiv 1\hspace{-6pt}
\pmod{N}}} \sum\limits_{\substack{d\mid x \\ d\mid c}} \mu (d) &
=\sum\limits_{\substack{ x\in I_1 \\ x\equiv 1\hspace{-6pt}
\pmod{N}}} \sum\limits_{\substack{d\mid x,\, d\mid c \\ (d,N)=1}}
\mu (d) =\sum\limits_{\substack{d\mid c \\ (d,N)=1}} \mu (d)
\sum\limits_{\substack{x\in I_1,\, d\mid x \\
x\equiv 1\hspace{-6pt} \pmod{N}}} 1 \\
& =\sum\limits_{\substack{d\mid c \\ (d,N)=1}} \mu (d) \left(
\frac{\vert I_1\vert}{dN}+O(1)\right) =\frac{\vert I_1\vert}{N}
\,\Phi_N(c) +O_\eps (c^\eps),
\end{split}
\end{equation*}
which completes the proof.
\end{proof}

Denote by $V_I(f)$ the total variation of a function $f$ defined
on the interval $I$.

\begin{cor}\label{C3}
For $I$ interval, $C^1$ functions $f_1,f_2:I\rightarrow \R$ with
$f_1\leq f_2$, and $T\geq 1$ integer, the cardinality
$\NN_{c,N}(f_1,f_2)$ of the set
\begin{equation*}
\big\{ (a,d)\in\Z^2 : d\in I,\ f_1(d)\leq a\leq f_2(d),\ a\equiv
1,\ d\equiv 1 \hspace{-6pt}\pmod N,\ ad\equiv
1\hspace{-6pt}\pmod{Nc} \big\}
\end{equation*}
can be expressed as
\begin{equation*}
\NN_{c,N} (f_1,f_2)=\frac{\Phi_N(c)}{\vert c\vert N^2} \int_I
(f_2-f_1) +\EE_{c,N,f_1,f_2}\qquad (\vert
c\vert\rightarrow\infty),
\end{equation*}
with
\begin{equation*}
\EE_{c,N,f_1,f_2} \ll_{\eps,N} \frac{\vert I\vert}{T\vert c\vert}
\big( V_I(f_1)+V_I(f_2)\big)+T\vert c\vert^{1/2+\eps} \left(
1+\frac{\vert I\vert}{T\vert c\vert}\right)\left( 1+\frac{\|
f_1\|_\infty+\| f_2\|_\infty}{\vert c\vert}\right).
\end{equation*}
\end{cor}

\begin{proof}
This follows from Proposition \ref{P1} as in the proof of
\cite[Lemma 3.1]{Bo}.
\end{proof}

\begin{lemma}\label{L2}
For every interval $J$ and every $C^1$ function $f:J\rightarrow
\R$
\begin{equation*}
\sum\limits_{c\in J} \Phi_N(c)f(c)=C_N \int_J f+O\Big(  \big( \|
f\|_\infty +V_J (f)\big)\log \sup\limits_{\xi\in
J}\vert\xi\vert\Big),
\end{equation*}
with $C_N$ as in \eqref{1.3}.
\end{lemma}

\begin{proof}
We can assume without loss of generality that $J=(0,Q]$. For each
$n\geq 1$ consider the $n$-dilate function $f_n (x):=f(nx)$, $x\in
[0,\frac{Q}{n}]$, for which $\| f_n\|_\infty =\| f\|_\infty$,
$\int_0^{Q/n} f_n=\int_0^Q f$, and $V_0^{Q/n}(f_n)=V_0^Q (f)$.
Using M\" obius and Euler-MacLaurin summation we get
\begin{equation*}
\begin{split}
\sum\limits_{c=1}^Q \Phi_N(c) f(c) & = \sum\limits_{c=1}^Q
\sum\limits_{\substack{n\mid c \\ (n,N)=1}} \frac{\mu (n)}{n} f(c)
=\sum\limits_{\substack{n\leq Q \\ (n,N)=1}}
\frac{\mu (n)}{n} \sum\limits_{c\leq Q/n} f_n (c) \\
& =\sum\limits_{\substack{n\leq Q \\ (n,N)=1}} \frac{\mu (n)}{n}
\left( \int_0^{Q/n} f_n +O\big( \|
f_n\|_\infty+V_0^{Q/n}(f)\big)\right) \\
& =\Bigg( \sum\limits_{\substack{n\geq 1 \\ (n,N)=1}} \frac{\mu
(n)}{n^2} +O\Big( \frac{1}{Q}\Big)\Bigg)\int_0^Q f+ O\Big(
\log Q \big( \| f\|_\infty+V_0^Q (f)\big)\Big) \\
& =C_N \int_0^Q f+O\Big( \log Q \big( \| f\|_\infty+V_0^Q
(f)\big)\Big),
\end{split}
\end{equation*}
which represents the desired conclusion.
\end{proof}

\begin{cor}\label{C4}
For every interval $I$ and every $C^1$ function $f:I\rightarrow
\R$
\begin{equation*}
\sum\limits_{\substack{c\in I \\ N\mid c}} \Phi_N
(c)f(c)=\frac{C_N}{N} \int_I f +O\Big(  \big( \| f\|_\infty +V_I
(f)\big)\log \sup\limits_{\xi\in I}\vert\xi\vert\Big).
\end{equation*}
\end{cor}

\begin{proof}
Apply Lemma \ref{L2} to $J=\frac{1}{N}I$, $f_N(x)=f(Nx)$, using
$\Phi_N (Nc^\prime)=\Phi_N(c^\prime)$, $\int_J f_N=\frac{1}{N}
\int_I f$, $\| f_N\|_\infty=\| f\|_\infty$, and $V_I (f_N)=V_J
(f)$.
\end{proof}

\section{Proof of the main results}
We first estimate the quantity defined in \eqref{2.10}.

\begin{proposition}\label{P2}
For every positive integer $N$ and every $z_0\in \H$, $\beta\in
[-\infty,\infty]$, as $Q\rightarrow \infty$,
\begin{equation*}
\mathfrak{N}_Q(N,z_0;\beta)=\frac{\pi(\pi+2\arctan\beta)C_N}{2N^3}Q^2
+O_{\eps,N,z_0}(Q^{7/4+\eps}).
\end{equation*}
\end{proposition}

\begin{proof}
Define
\begin{equation*}
\begin{split}
& I_c =cx_0+\begin{cases} \big[ -\sqrt{Q^2-c^2y_0^2},\ \min\{
\beta cy_0,\sqrt{Q^2-c^2 y_0^2}\ \}\big] & \mbox{\rm if
$c\in \left[ 0,Q/y_0 \right],$}\\
\big[ \max \{ \beta cy_0,-\sqrt{Q^2-c^2y_0^2}\ \}, \sqrt{Q^2-c^2
y_0^2}\ \big] & \mbox{\rm if $c\in \left[-Q/y_0,0\right],$}
\end{cases} \\
& f(c,a)=\vert c\vert y_0\sqrt{\frac{Q^2}{c^2y_0^2+(a-cx_0)^2}-1},
\\ & f_1(c,a)=-cx_0-f(c,a),\quad f_2(c,a)=-cx_0+f(c,a), \qquad c\in
\left[-Q/y_0,Q/y_0 \right],\ a\in I_c,\\
& F(c)=F_{z_0,\beta}(c)=\frac{2}{\vert c\vert} \int_{I_c} f(c,a)\
da.
\end{split}
\end{equation*}

Writing the inequalities from \eqref{2.10} as
\begin{equation*}
\begin{cases}
\vert C\vert \leq Q,\\
-\sqrt{Q^2-C^2}\leq A\leq \sqrt{Q^2-C^2}\quad \mbox{\rm and} \quad
\begin{cases} A\leq \beta C & \mbox{\rm if $C>0$,} \\ A\geq \beta
C & \mbox{\rm if $C<0$,}
\end{cases} \\
-\vert C\vert \sqrt{\frac{Q^2}{C^2+A^2}-1} \leq D\leq \vert C\vert
\sqrt{\frac{Q^2}{C^2+A^2}-1},
\end{cases}
\end{equation*}
and using \eqref{2.1} we gather
\begin{equation}\label{4.1}
\begin{split}
\mathfrak{N}_Q(N,z_0;\beta) & =\# \Big\{ \gamma\in\Gamma(N): \vert
c\vert y_0\leq Q,\ a\in I_c,\ d\in
[ f_1(c,a),f_2(c,a)]\Big\} \\
& =\sum\limits_{\vert c\vert\leq Q/y_0} \NN_{c,N}\big(
f_1(c,\cdot),f_2 (c,\cdot)\big).
\end{split}
\end{equation}

Note that $\max\big\{ \|
f(c,\cdot)\|_\infty,V_{I_c}\big(f(c,\cdot)\big)\big\}\ll Q$ on
$I_c$, thus Corollary \ref{C3} with $T=[Q^{1/4}]$ gives
\begin{equation}\label{4.2}
\NN_{c,N}\big( f_1 (c,\cdot),f_2(c,\cdot)\big)
=\frac{1}{N^2}\Phi_N(c)F(c)+\EE_{c,N},
\end{equation}
with
\begin{equation}\label{4.3}
\EE_{c,N} \ll_{\eps,N} Q^{7/4}\vert c\vert^{-1}+Q^{5/4} \vert
c\vert^{-1/2+\eps} +Q^2 \vert c\vert^{-3/2+\eps}.
\end{equation}

Fix some constant $\alpha\in \big[\frac{1}{2},\frac{3}{4}\big]$.
The relation $bc\equiv -1 \pmod{\vert a\vert}$ and the constraint
$\vert a\vert \ll_{z_0} Q$ give the trivial estimate
\begin{equation}\label{4.4}
\sum\limits_{\vert c\vert\leq Q^\alpha} \NN_{c,N}\big(
f_1(c,\cdot),f_2(c,\cdot)\big) \ll_{z_0} \sum\limits_{1\leq\vert
a\vert \leq Q} Q^\alpha\frac{Q}{\vert a\vert}\ll Q^{1+\alpha}\log
Q \ll_\eps Q^{7/4+\eps}.
\end{equation}
On the other hand \eqref{4.3} leads to
\begin{equation}\label{4.5}
\begin{split}
\sum\limits_{Q^\alpha<\vert c\vert\leq Q/y_0} \EE_{c,N} &
\ll_{\eps,z_0,N} Q^{7/4}\log Q+Q^{5/4}\sum\limits_{1\leq c\leq Q}
c^{-1/2+\eps}+Q^2\sum\limits_{c>Q^\alpha} c^{-3/2+\eps} \\
& \ll_{\eps}
Q^{7/4+\eps}+Q^{5/4+1/2+\eps}+Q^{2+\alpha(-1/2+\eps)}\ll
Q^{7/4+\eps}.
\end{split}
\end{equation}

From \eqref{4.1}--\eqref{4.5} we now infer
\begin{equation}\label{4.6}
\mathfrak{N}_Q(N,z_0;\beta) =\frac{1}{N^2} \sum\limits_{Q^\alpha
\leq\vert c\vert\leq Q/y_0} \Phi_N(c)F(c)+O_{\eps,N,z_0}
(Q^{7/4+\eps}).
\end{equation}

Using $I_c \subseteq
\big[-\sqrt{Q^2-c^2y_0^2},\sqrt{Q^2-c^2y_0^2}\,\big]$ and the
change of variable $u=C\tan x$ we get
\begin{equation*}
\begin{split}
F(c) & =2y_0\int_{I_c} \sqrt{\frac{Q^2}{c^2y_0^2+(a-cx_0)^2}-1}\,
da \leq 4y_0 \int_0^{\sqrt{Q^2-C^2}} \sqrt{\frac{Q^2}{C^2+u^2}-1}\, du \\
& =4y_0 \int_0^{\arctan \sqrt{Q^2/C^2-1}}
\sqrt{Q^2-\frac{C^2}{\cos^2 x}}\, \frac{dx}{\cos x}\leq 4y_0 Q
\int_0^{\arctan \sqrt{Q^2/C^2-1}} \frac{dx}{\cos x} \\
& =2y_0Q\log \frac{1+\sin x}{1-\sin x}\
\bigg|_{x=0}^{\arctan\sqrt{Q^2/C^2 -1}} =4y_0 Q \log \bigg(
\frac{Q}{C}+\sqrt{\frac{Q^2}{C^2}-1}\ \bigg)\ll_{z_0} Q\log Q.
\end{split}
\end{equation*}
The total variation of $F$ on $\big[-\frac{Q}{y_0},-Q^\alpha\big]$
and on $\big[ Q^\alpha,\frac{Q}{y_0}\big]$ is also $\ll_{z_0}
Q\log Q$ because $F$ is slowly oscillating. Applying Corollary
\ref{C4} to the sum from \eqref{4.6} we now infer
\begin{equation*}
\begin{split}
\mathfrak{N}_Q(N,z_0;\beta) & =\frac{C_N}{N^3} \int_{Q^\alpha \leq
\vert c\vert\leq Q/y_0} F(c)\, dc +O_{\eps,N,z_0}(Q^{7/4+\eps}) \\
& =\frac{C_N}{N^3} \int_{-Q/y_0}^{Q/y_0} F(c)\,
dc+O_{\eps,N,z_0}(Q^{7/4+\eps}).
\end{split}
\end{equation*}
Using the substitution $c=\frac{Qu}{y_0}$,
$a=Qv+cx_0=\big(v+\frac{ux_0}{y_0}\big) Q$, the integral in the
main term above is evaluated as
\begin{equation*}
\begin{split}
\int_{-Q/y_0}^{Q/y_0} F(c)\, dc & =2\int_{-Q/y_0}^{Q/y_0}
\int_{I_c} f(c,a)\, da\, dc
\\ & =2\iint_{\substack{u^2+v^2\leq 1\\ u\geq 0,\, v\leq\beta u}}
\sqrt{\frac{1}{u^2+v^2}-1} \ du\, dv+
2\iint_{\substack{u^2+v^2\leq 1\\ u\leq 0,\, v\geq\beta u}}
\sqrt{\frac{1}{u^2+v^2}-1}
\ du\, dv\\
& =2\int_0^1 \int_{-\pi/2}^{\arctan\beta} \sqrt{1-r^2} \
d\theta\,dr+2\int_0^1\int_{\pi/2}^{\pi+\arctan\beta} \sqrt{1-r^2}\
d\theta\,dr \\
& =\frac{\pi(\pi+2\arctan\beta)}{2}.
\end{split}
\end{equation*}
This completes the proof of the proposition.
\end{proof}

Taking stock on \eqref{2.9} we obtain (recall that
$Q^2=e^R+O(e^{-R})$)
\begin{equation}\label{4.7}
\begin{split}
\# \Gamma(N)_{z_0,R} & =\# \left\{ \gamma\in\Gamma(N):
A^2+B^2+C^2+D^2\leq Q^2\right\}
=\mathfrak{N}_{\sqrt{Q^2+O_{z_0}(Q)}} (N,z_0;\infty) \\
& =\frac{\pi^2C_NQ^2}{N^3}+O_{\eps,N,z_0}(Q^{7/4+\eps})=
\frac{6Q^2}{[\Gamma(1):\Gamma(N)]} +O_{\eps,N,z_0}(Q^{7/4+\eps})\\
& =\frac{6e^R}{[\Gamma(1):\Gamma(N)]}+O_{\eps,N,z_0}\left(
e^{(7/8+\eps)R}\right),
\end{split}
\end{equation}
which proves Corollary \ref{C2}.

\begin{proof}[Proof of Theorem \ref{T1}]
Set $\mathfrak{N}_Q (\beta)=\mathfrak{N}_Q (N,z_0;\beta)$. As a
consequence of Proposition \ref{P2} and of the inequality $\vert
\arctan (\beta+\beta_0)-\arctan \beta \vert \leq \vert
\beta_0\vert$ we have
\begin{equation}\label{4.8}
\left| \mathfrak{N}_Q (\beta+\beta_0)-\mathfrak{N}_Q
(\beta)\right| \ll_{\eps,N,z_0} Q^2 \vert\beta_0\vert
+Q^{7/4+\eps}.
\end{equation}

Let $S_Q(\omega)=S_Q (N,z_0,z_1;\omega)$ denote the cardinality of
the set of $\gamma\in\Gamma (N)$ with $A^2+B^2+C^2+D^2 \leq Q^2$
and $-\frac{\pi}{2}\leq \theta_{z_0,z_1}(\gamma)\leq \omega$.
Partitioning this set according to whether or not $\min \{ \vert
A\vert,\vert C\vert\} > Q^\alpha$ and employing Lemma \ref{L1} we
find that, up to an error $\ll_{z_0} Q^{1+\alpha}\log Q$,
$S_Q(\omega)$ equals
\begin{equation}\label{4.9}
\# \left\{ \gamma\in\Gamma (N): A^2+B^2+C^2+D^2 \leq Q^2, \vert
A\vert,\vert C\vert >Q^\alpha, -\pi/2\leq\theta_{z_0,z_1}
(\gamma)\leq \omega \right\}.
\end{equation}
By \eqref{2.9} there is $K_2=K_2(z_0)>0$ such that the number in
\eqref{4.9} is
\begin{equation}\label{4.10}
\leq \#\left\{ \gamma\in\Gamma (N): (C^2+A^2)\bigg(
1+\frac{D^2}{C^2}\bigg)\leq Q_1^2, \vert A\vert ,\vert C\vert
>Q^\alpha,\ -\frac{\pi}{2}\leq \theta_{z_0,z_1}(\gamma)\leq \omega
\right\} ,
\end{equation}
where we set $Q_1:=\sqrt{Q^2+K_2 Q}=Q+O_{z_0}(1)$. According to
\eqref{2.8} the number in \eqref{4.10} is
\begin{equation*}
\leq \#\left\{ \gamma\in\Gamma (N): (C^2+A^2)\bigg(
1+\frac{D^2}{C^2} \bigg) \leq Q_1^2 ,\ \frac{A}{C}\in x_*
+\mathfrak{S} (y_*,\tan\omega)
+\bigg[-\frac{3K_1}{Q^{2\alpha}},\frac{3K_1}{Q^{2\alpha}}\bigg]\right\}.
\end{equation*}
Taking $\alpha=\frac{1}{8}$ and applying \eqref{4.8} to
$\vert\beta_0\vert =Q^{-2\alpha}=Q^{-1/4}$ we find
\begin{equation}\label{4.11}
\begin{split}
S_Q(\omega)\leq & \# \left\{ \gamma\in\Gamma (N): (C^2+A^2)\bigg(
1+\frac{D^2}{C^2}\bigg)\leq Q_1^2 ,\ \frac{A}{C}\in
x_*+\mathfrak{S} (y_*,\tan\omega)\right\} \\ & \qquad
+O_{\eps,N,z_0,z_1} (Q^{7/4+\eps}).
\end{split}
\end{equation}
The number of matrices $\gamma\in\Gamma (N)$ for which
$\frac{A}{C}=\mu$ and $A^2+B^2+C^2+D^2\leq Q^2$ is $\ll_{z_0,\mu}
Q$ as $Q\rightarrow \infty$. Using this fact together with
\eqref{2.7}, \eqref{2.9}, and \eqref{2.10}, we find that, up to a
term of order $O_{z_0}(Q_1)=O_{z_0}(Q)$, the main term in the
right-hand side of \eqref{4.11} is given by
\begin{equation}\label{4.12}
\begin{split}
& \begin{cases} \mathfrak{N}_{Q_1} \big( x_*-y_* \cot
\frac{\omega}{2}\big)+\mathfrak{N}_{Q_1} \big( x_*+y_*
\tan\frac{\omega}{2}\big)-\mathfrak{N}_{Q_1} (x_*-y_*) & \\
\qquad +\mathfrak{N}_{Q_1}
(\infty)-\mathfrak{N}_{Q_1} (x_*+y_*) & \mbox{\rm if $\omega >0$}, \\
\mathfrak{N}_{Q_1} (x_*)-\mathfrak{N}_{Q_1}
(x_*-y_*)+\mathfrak{N}_{Q_1}
(\infty)-\mathfrak{N}_{Q_1} (x_*+y_*) & \mbox{\rm if $\omega =0$}, \\
\mathfrak{N}_{Q_1} \big( x_*+y_*\tan\frac{\omega}{2}\big)
-\mathfrak{N}_{Q_1} (x_*-y_*) \\ \qquad +\mathfrak{N}_{Q_1} \big(
x_*-y_* \cot\frac{\omega}{2}\big) -\mathfrak{N}_{Q_1} (x_*+y_*) &
\mbox{\rm if $\omega <0$},
\end{cases} \\
& \qquad =\mathfrak{N}_{Q_1} \left( x_*+y_*
\tan\frac{\omega}{2}\right) +\mathfrak{N}_{Q_1} \left( x_*-y_*
\cot\frac{\omega}{2}\right) -\mathfrak{N}_{Q_1} (x_*+y_*)
-\mathfrak{N}_{Q_1} (x_*-y_*) \\ & \qquad \qquad +\begin{cases}
\mathfrak{N}_{Q_1} (\infty) & \mbox{\rm if $\omega >0$}, \\
0 & \mbox{\rm if $\omega <0$}. \end{cases}
\end{split}
\end{equation}
As a result of Proposition \ref{P2} and $Q_1=Q+O_{z_0}(1)$ the
expression in \eqref{4.12} equals
\begin{equation*}
\begin{split}
& \frac{\pi C_N Q^2}{N^3}\Bigg( \arctan \left(
x_*+y_*\tan\frac{\omega}{2}\right) +\arctan \left( x_* -y_* \cot
\frac{\omega}{2} \right) -\arctan (x_*+y_*) \\ & \hspace{4cm}
-\arctan (x_* -y_*)+\begin{cases} \pi & \mbox{\rm if $\omega
>0$,} \\ 0 & \mbox{\rm if $\omega <0$,} \end{cases}
\Bigg)+O_{\eps,N,z_0,z_1}(Q^{7/4+\eps}) .
\end{split}
\end{equation*}
Letting $\Xi_{x_*,y_*}$ as in \eqref{1.4} we now infer
\begin{equation}\label{4.13}
S_Q(\omega) \leq \frac{\pi^2 C_N \Xi_{x_*,y_*} (\omega)}{N^3}
Q^2+O_{\eps,N,z_0,z_1} (Q^{7/4+\eps}).
\end{equation}
The opposite inequality
\begin{equation*}
S_Q(\omega) \geq \frac{\pi^2 C_N \Xi_{x_*,y_*} (\omega)}{N^3}
Q^2+O_{\eps,N,z_0,z_1} (Q^{7/4+\eps})
\end{equation*}
is derived in a similar way. Therefore equality holds in
\eqref{4.13}. Equality \eqref{1.5} now follows taking $Q^2=2\cosh
R=e^R+e^{-R}$.

Estimates \eqref{1.5} and \eqref{4.7} provide
\begin{equation}\label{4.14}
\begin{split}
\P_{\Gamma (N), B(z_0,R),z_0,z_1} (\omega) & = \frac{\#\big\{
\gamma\in \Gamma (N)_{z_0,R}:
-\pi/2\leq \theta_{z_0,z_1} (\gamma) \leq\omega\big\}}{\# \Gamma (N)_{z_0,R}}\\
& =\frac{\frac{\pi^2 C_N}{N^3}
\Xi_{x_*,y_*}(\omega)e^R+O_{\eps,N,z_0,z_1}(e^{(7/8+\eps)R})}{\frac{\pi^2
C_N}{N^3} e^R+O_{\eps,N,z_0,z_1}(e^{(7/8+\eps)R})}
\\ & =\Xi_{x_*,y_*}(\omega)+O_{\eps,N,z_0,z_1} \left(
e^{(-1/8+\eps)R}\right).
\end{split}
\end{equation}

The function $\Xi_{x_*,y_*}$ is differentiable on
$[-\frac{\pi}{2},\frac{\pi}{2}]$ with
\begin{equation}\label{4.15}
\begin{split}
\Xi_{x_*,y_*}^\prime (\omega) &
=\frac{y_*}{2\pi\cos^2\frac{\omega}{2}\big(1+(x_*+y_*\tan\frac{\omega}{2})^2\big)}
+\frac{y_*}{2\pi\sin^2\frac{\omega}{2}\big(
1+(x_*-y_*\cot\frac{\omega}{2})^2\big)}\\
& =\frac{y_*}{2\pi}\left( \frac{1}{\cos^2\frac{\omega}{2}+\big(
x_*\cos\frac{\omega}{2}+y_*\sin\frac{\omega}{2}\big)^2}+
\frac{1}{\sin^2\frac{\omega}{2}+\big(x_*\sin\frac{\omega}{2}-y_*\cos\frac{\omega}{2}\big)^2}
\right)\\
& =\frac{2}{\pi}\cdot\frac{y_*(1+x_*^2+y_*^2)}{(1+x_*^2+y_*^2)^2
-\big( (1+x_*^2-y_*^2)\cos\omega +2x_*y_*\sin\omega\big)^2} \\
& =\frac{1}{\pi} \varrho_{z_0,z_1} (\omega).
\end{split}
\end{equation}
The second part of Theorem \ref{T1} now follows from \eqref{4.14}
and \eqref{4.15}.
\end{proof}

\section*{Acknowledgments}
I am grateful to Professor Peter Sarnak for his interest and for
bringing to my attention references \cite{Go} and \cite{Sha}, and
to Professor Freydoon Shahidi for kindly providing me a reprint of
\cite{Sha}.

\end{document}